\documentclass[12pt,thmsa]{article}
\usepackage{amsmath}
\usepackage{graphicx}
\usepackage{amsfonts}
\usepackage{amssymb}

\newcommand{\R}{\mathbb{R}}
\newcommand{\Z}{\mathbb{Z}}
\newcommand{\N}{\mathbb{N}}
\newcommand{\p}{\mathbb{P}}
\newcommand{\E}{\mathbb{E}}

\newtheorem{thm}{Theorem}
\newtheorem{assumption}{Assumption}
\newtheorem{cor}{Corollary}
\newtheorem{lem}{Lemma}

\begin{document}

\begin{center}
{\Large \textbf{On estimation and prediction in  spatial functional linear regression model}}

\bigskip

St\'ephane BOUKA$^1$, Sophie DABO-NIANG$^{2,3}$   and Guy Martial  NKIET$^1$

\bigskip

$^1$ URMI, Universit\'{e} des Sciences et Techniques de Masuku, Franceville, Gabon.

$^2$Laboratoire LEM, CNRS 9221, Universit\'e de  Lille, France.

$^3$INRIA-MODAL, Lille, France. 

\bigskip

E-mail : stephane.bouka@mathsinfo.univ-masukul.com;  sophie.dabo@univ-lille3.fr; guymartial.nkiet@mathsinfo.univ-masuku.com.

\bigskip
\end{center}

\noindent\textbf{Abstract.}We consider a spatial functional linear regression, where a scalar response is related to a square integrable spatial functional process. We use a smoothing spline estimator for the functional slope parameter and establish a finite sample bound for variance of this estimator under mixing spatial dependence. Then, we give a bound of the prediction error. Finally, we illustrate our results by simulations.

\bigskip



\section{Introduction}

\noindent Consider the following spatial functional linear regression model where the spatial scalar response $(Y_{\mathbf{i}}\in\R,\ \mathbf{i}\in D\subset\Z^{d})$ is related to a square integrable spatial functional process $(X_{\mathbf{i}}\in{\mathcal F},\ \mathbf{i}\in D\subset\Z^{d})$ through
\begin{equation}\label{rl1.1}
Y_{\mathbf{i}}=\beta_{0}+\int_{I}\beta(t)X_{\mathbf{i}}(t)dt+\epsilon_{\mathbf{i}}\ ,\ \ \mathbf{i}\in\Z^{d}
\end{equation}
where $\beta_{0}$ is a constant, $I$ is the domain of $X_{\mathbf{i}}$, ${\mathcal F}$ is a space of functions endowed with a semi-norm, $\beta$ is an unknown function representing the slope function, and $\left(\epsilon_{\mathbf{i}}\right)_{\mathbf{i}\in\Z^{d}}$ is a centered random spatial noise and with variance $\sigma^{2}_{\epsilon}>0$.
The functional linear regression with functional or scalar response has been the focus of various investigations. There exist many contributions in this field for non spatial data, and recent references are: \cite{aue}, \cite{ComteandJohannes12}, \cite{Crambesetal09}, \cite{Cuevas14}, \cite{hormann_kokoszka10}, \cite{kneipal16}, \cite{liuetal17}, \cite{mas_pumo09}, \cite{zhouetal16}. This work is motivated by a large number of applications for which the data are of spatial nature. For example, non-parametric prediction from kriging methods for geostatistical functional data was tackled in \cite{bohorquezetal16},  \cite{bohorquezetal17}, \cite{Giraldoetal12}, \cite{giraldoetal14}, \cite{giraldoetal18}, \cite{giraldoetal11} and \cite{nerinietal10} whereas spatial autoregressive functional  models were considered in \cite{ruiz2011spatial, ENV:ENV1143}. In this paper, we are interested in estimation of the slope function $\beta$ in model (\ref{rl1.1}). To the best of our knowledge, this problem have not yet been considered for the basic spatial functional linear regression model, but only for non spatial data (e.g  \cite{Crambesetal09}) or for spatial linear regression model with derivatives (see \cite{boukaetal18}). The paper is organized as follow. The section 2 is devoted to the estimation of the estimator that will be use. Assumptions and main results are stated in Section 3, and a simulation study is given in Section 4. The proofs are postponed to Section 5.

\section{Smoothing splines estimation of slope function}\label{section r2}

\noindent In this section, we give an estimator of $\beta$ in (\ref{rl1.1}) by using an approach similar to the one of  \cite{Crambesetal09}. Since this procedure of estimation does not take into account the nature of the dependence of the data, we obtain an  estimator that has the same form than that  of   \cite{Crambesetal09}. The process $(X_{\mathbf{i}}, Y_{\mathbf{i}})_{\mathbf{i}\in\Z^d}$ is defined on probability space $(\Omega,{\mathcal A}, \p)$ with the same distribution as a couple of variable $(X,Y)$. For $\mathbf{n}=(n,\cdots,n)$ with $n\in\N^*$, let ${\mathcal I}_{\mathbf{n}}:=\{1,\cdots,n\}^d$ be a grid of points in $\Z^d$ and consider observations $(X_{\mathbf{i}}, Y_{\mathbf{i}})_{\mathbf{i}\in{\mathcal I}_{\mathbf{n}}}$. We assume that the random functions $X_{\mathbf{i}}$ are observed at $p$ equidistant points $t_{1},...,t_{p}\in I:=[0,1]$, with $t_{j}=\frac{j}{p}$ for all $j=1,...,p$. By using the lexico-graphic order, the previous sample is rewritten as $\{(X_{\mathbf{i}_i},Y_{\mathbf{i}_i})\}_{1\le i\le n^d}$,  then we put $\mathbf{Y}=(Y_{\mathbf{i}_{1}}-\overline{Y},...,Y_{\mathbf{i}_{n^d}}-\overline{Y})^T$ (where $u^T$ denotes the transposed of $u$) and we consider the $n^d\times p$ matrix $\mathbf{X}$ with general term $X_{\mathbf{i}_{i}}(t_{j})-\overline{X}(t_{j})$ for $i=1,...,n^d$, $j=1,...,p$. Then, we consider the estimator $\widehat{\beta}$ of $\beta$ given by
\begin{eqnarray}
\widehat{\beta}(t)&=&\mathbf {D}(t)^{T}(\mathbf{D}^T\mathbf{D})^{-1}\mathbf{D}^T\widehat{\boldsymbol{\beta}}\label{sfr1}\\
 \text{with}\ \ \ \widehat{\boldsymbol{\beta}}&=&\frac{1}{n^d}\left(\frac{1}{n^dp}\mathbf{X}^T\mathbf{X}+\rho\mathbf{A}_{m}\right)^{-1}\mathbf{X}^T\mathbf{Y} ,\label{r2.2}
\end{eqnarray}
where $\rho>0$ is a smoothing parameter, $\mathbf{A}_{m}$ is a $p\times p$ symmetric matrix defined from B-splines (see \cite{Crambesetal09} for details), $\mathbf{D}(t)=(D_1(t),\cdots,D_p(t))^T$ is a functional basis of the $p$-dimensional linear space $NS^m(t_1,\cdots,t_p)$  of functions $v$ having a $m$-th order derivative  $v^{(m)}$ that belongs to $L^2([0,1])$, and $\mathbf{D}$ is the $p\times p$ matrix with general term $D_i(t_j)$ for  $i,j=1,\cdots,p$. For estimating the intercept $\beta_{0}$ we take $\widehat{\beta_{0}}=\overline{Y}-\left\langle\widehat{\beta},\overline{X}\right\rangle$, where $\left\langle.,.\right\rangle$ denotes the usual inner product of $L^2([0,1])$.

\section{Assumptions and main results}\label{section r3}

\noindent In this section, we first introduce the assumptions that are needed to obtain the main  results of the paper, then theorems that give the rate of convergence of the estimator $\widehat{\beta}$ and also that of  the prediction at a non-visited site are established.

\subsection{Assumptions}

\begin{assumption}\label{ar2}
$\beta$ is $m$-times differentiable and $\beta^{(m)}$ belongs to $L^{2}([0,1])$.
\end{assumption}

\begin{assumption}\label{ar11} There exist $\kappa\in ]0,1[$, $\delta_1>0$ and $C_1>0$ such that,  
for any $(t,s)\in I^2$,  
$\p\left(|X(t)-X(s)|\leq C_{1}|t-s|^{\kappa}\right)\geq1-\delta_{1}.
$
\end{assumption}
\begin{assumption}\label{ar12}
  For $C_{2}\in\mathbb{R}_+^\ast$  and all $r\in\mathbb{N}^\ast$ there exists a $r$-dimensional linear subspace ${\mathcal L}_{r}$ of $L^{2}([0,1])$ and a real $q\in ]0,1[$ such that
\[
\E\left(\inf_{f\in{\mathcal L}_{r}}\sup_{t}|X(t)-f(t)|^{2} \right)\leq C_{2}r^{-2q}.
\]
\end{assumption}

\begin{assumption}\label{ar13}
For any $(j,\ell)\in\mathbb{N}^\ast$,
\begin{eqnarray*}
Var\left(\frac{1}{n^d}\sum^{n^d}_{i=1}\left\langle X_{\mathbf{i}_{i}}-\E(X),\zeta_{j}\right\rangle\left\langle X_{\mathbf{i}_{i}}-\E(X),\zeta_{\ell}\right\rangle\right)
\\\leq\frac{C_{3}}{n^d}\E\left(\left\langle X-\E(X),\zeta_{j}\right\rangle^{2}\right)\E\left(\left\langle X-\E(X),\zeta_{\ell}\right\rangle^{2}\right)
\end{eqnarray*}
where $0<C_{3}<\infty$ and  $\left\{\zeta_{j}\right\}_{j\in\mathbb{N}^\ast}$ is a complete orthonormal system of eigenfunctions of the operator $\Gamma$ from $L^2([0,1])$ to itself defined by:
\[
\Gamma u:=\E\left(<u,X-\E(X)>(X-\E(X))\right),
\] 
each $\zeta_j$ being associated with the $j$-th largest eigenvalue $\lambda_j$.
\end{assumption}

\noindent Assumptions \ref{ar2}--\ref{ar13} are technical conditions that are similar to the ones considered in \cite{Crambesetal09}. In order  to give the remaining assumptions, let us first recall the notion of polynomial mixing dependence. Letting $\alpha$ be the $\alpha$-mixing coefficient given, for two sub $\sigma$-algebras ${\mathcal U}$ and ${\mathcal V}$ of ${\mathcal A}$, by
\begin{eqnarray*}
\alpha({\mathcal U},{\mathcal V})=\sup\{|\p(A\cap B)-\p(A)\p(B)|, A\in{\mathcal U},B\in{\mathcal V}\},
\end{eqnarray*}
we  consider the strong mixing coefficient (see \cite{mach1}) related to a random field $(Z_{\mathbf{i}})_{\mathbf{i}\in\Z^d}$, defined as
\begin{eqnarray}\label{ar3.5}
\alpha_{1,\infty}(u)=\sup\{\alpha(\sigma(Z_{\mathbf{i}}),F_{\Lambda}), \mathbf{i}\in\Z^{d},\Lambda\subset\Z^{d}, \delta(\Lambda,\{\mathbf{i}\})\geq u\},
\end{eqnarray}
where  $ F_{\Lambda}=\sigma(Z_{\mathbf{i}};\mathbf{i}\in{\Lambda})$ and the distance $\delta$ is defined for any subsets $\Gamma_{1}$ and $\Gamma_{2}$ of $\Z^{d}$ by $\delta(\Gamma_{1},\Gamma_{2})=\min\{||\mathbf{i}-\mathbf{j}||_{2}, \mathbf{i}\in\Gamma_{1}, \mathbf{j}\in\Gamma_{2}\}$ where $\|.\|_{2}$ is the usual Euclidean norm of $\R^d$. Then, $(Z_{\mathbf{i}})_{\mathbf{i}\in\Z^d}$ is polynomial mixing  if the related strong mixing coefficients satisfy $\alpha_{1,\infty}(u)=O(u^{-\theta})$, $\theta>0$. 
\begin{assumption}\label{ar22}
$\{\epsilon_{\mathbf{i}}\}_{\mathbf{i}\in\Z^{d}}$ is a strictly stationary random field, polynomial mixing, independent of $\{X_{\mathbf{i}}\}_{\mathbf{i}\in\Z^{d}}$ and  such that $\sup_{\mathbf{i}\in\Z^{d}}\left|\epsilon_{\mathbf{i}} \right|<M_1$ almost surely,  where $M_1$ is a strictly positive constant.
\end{assumption}

\begin{assumption}\label{ar3}
$\left\{(X_{\mathbf{i}},Y_{\mathbf{i}})\right\}_{\mathbf{i}\in\Z^{d}}$ is a strictly stationary and polynomial mixing random field.
\end{assumption}

\begin{assumption}\label{ar4} There exists $M_2>0$ such that for all $\mathbf{i}\in\Z^{d}$, 
$\left\|X_{\mathbf{i}} \right\|<M_2$ almost surely.
\end{assumption}
 Assumptions \ref{ar22} and \ref{ar3} are classical  assumptions (see \cite{biau_cadre}). Assumption \ref{ar4} has already been made in some works (see, e.g., \cite{mas_pumo09}).
\begin{assumption}\label{9}
$X $ is an isotropic process such that for all $t$, $u$ in $[0,1]$,
$$Cov(X_{\mathbf{i}_i}(t),X_{\mathbf{i}_j}(u))=g(|t-u|)\ \varPsi(\delta(\{\mathbf{i}_i\},\{\mathbf{i}_j\}))\ \ \text{and}\ \ \varPsi(0)=1$$
where $g$ is a positive function and $\varPsi$ is a known $\R_+$-valued decreasing function that verified $\sum^\infty_{t=1}t^{d-1}\varPsi(t)<\infty$.
\end{assumption}
 The separable covariance structure stated in Assumption \ref{9} has also been used in \cite{liuetal17}. Examples on isotropic spatial models can be founded in  \cite{francisco}. We may mention for instance, the exponential spatial model.  

\subsection{The  results}\label{section r4}

We consider  the semi-norm $\|.\|_{\Gamma}$ defined by
\begin{eqnarray}\label{s1}
||u||^{2}_{\Gamma}:=\left\langle\Gamma u,u\right\rangle,\ \ u\in L^{2}([0,1]).
\end{eqnarray}
and the discretized empirical semi-norm defined for any $\mathbf{u}\in\R^{p}$ as
\begin{eqnarray*}
\|\mathbf{u}\|^{2}_{\Gamma_{n,p}}:=\frac{1}{p}\mathbf{u}^T\left(\frac{1}{n^dp}\mathbf{X}^T\mathbf{X}\right)\mathbf{u}.
\end{eqnarray*}
The following theorem gives   a bound of the estimator's variance. In this theorem, $\E_{\epsilon}$  refers  to the conditional expectation given $X_{\mathbf{i}_{1}},...,X_{\mathbf{i}_{n^d}}$.

\begin{thm}\label{tr1}
Under Assumptions \ref{ar2},  \ref{ar22}, \ref{ar3}  and \ref{ar4}  with $\alpha_{1,\infty}(u)=O(u^{-\theta})$, $\theta>d$, for all $\rho>n^{-2md}$ , if the eigenvalues $\lambda_{x,1}\ge\lambda_{x,2}\ge...\ge\lambda_{x,p}\ge0$ of $1/(n^dp)\mathbf{X}^T\mathbf{X}$ satisfy $\sum^{p}_{j=r+1}\lambda_{x,j}\leq C.r^{-2q}$ with $C>0$, $q>0$ and $r:=\lfloor\rho^{-1/(2m+2q+1)}\rfloor$, then
\begin{eqnarray}
\E_{\epsilon}(\|\widehat{\boldsymbol{\beta}}-\E_{\epsilon}(\widehat{\boldsymbol{\beta}})\|^{2}_{\Gamma_{n,p}})\leq
\left(\frac{\sigma^{2}_{\epsilon}}{n^d}+\frac{c\ln n}{n^d}\right)\left(m+\lfloor\rho^{-1/(2m+2q+1)}\rfloor(2+C.C_{0})\right)\label{r3.2}
\end{eqnarray}
where $C_{0}>0$, $c>0$ and $\lfloor x\rfloor$ stands the integer part of $x$.
\end{thm}
\noindent Using Theorem $\ref{tr1}$ and Arguing as in \cite{Crambesetal09} , we obtain the Corollary below.
 
 \begin{cor}\label{cr1}
 Under assumptions of Theorem $\ref{tr1}$ together with \\Assumptions $\ref{ar11}$-$\ref{ar13}$, as well as $n^dp^{-2\kappa}=O(1)$, $\rho\rightarrow0$, $1/(n^d\rho)\rightarrow0$ as $n, p\rightarrow\infty$ we have
 \begin{gather}\label{r3.4}
 \|\widehat{\beta}-\beta\|^{2}_{\Gamma}=O_{p}\left(\rho +\left(n^d\rho^{1/(2m+2q+1)}\right)^{-1}\ln n+n^{-d(2q+1)/2}\right)
 \end{gather}
 \end{cor}
 \noindent Next, we give a bound for prediction error. For that, we assume what follows
 \begin{assumption}\label{8} The non-visited site $\mathbf{i}_0$ is such that
 $$\delta(\{\mathbf{i}_0\},\{\mathbf{i}_1,...,\mathbf{i}_n\})\ge \lfloor n^{2d/\theta}\rfloor$$
 \end{assumption}
\noindent In this Assumption \ref{8}, it is sufficient to choice $\theta$ large for doing the prediction at any non-visited site. 

\noindent we consider the  prediction $\widehat{Y}_{\mathbf{i}_{0}}$ and the ''theoretical'' prediction $Y^{*}_{\mathbf{i}_{0}}$ at a non-visited site $\mathbf{i}_{0}\in\Z^d$ such that $(X_{\mathbf{i}_0}, Y_{\mathbf{i}_0})$ has the same distribution than $(X, Y)$. In fact,
 \begin{eqnarray}\label{s2}
 \widehat{Y}_{\mathbf{i}_{0}}=\widehat{\beta}_{0}+\left\langle\widehat{\beta},X_{\mathbf{i}_{0}}\right\rangle\ \ and\ \ Y^{*}_{\mathbf{i}_{0}}=\beta_{0}+\left\langle\beta,X_{\mathbf{i}_{0}}\right\rangle
 \end{eqnarray}
 We are interested by the bound of the prediction error between $\widehat{Y}_{\mathbf{i}_{0}}$ and $Y^{*}_{\mathbf{i}_{0}}$.
 \begin{thm}\label{cr2}
 Suppose that assumptions of Corollary $\ref{cr1}$ together with assumptions $\ref{9}$--$\ref{8}$ hold. If $\sum_{j\ge1}\lambda^{1/4}_j<\infty$, $2q>1$, $\rho\sim n^{-d(2m+2q+1)/(2m+2q+2)}$ and $p$ is chosen sufficiently large compared to $n^d$, then 
 \begin{eqnarray*}
 &&\E\left((\widehat{Y}_{\mathbf{i}_{0}}-Y^{*}_{\mathbf{i}_{0}})^{2}|\widehat{\beta}_{0},\widehat{\beta}\right)=
 O_{p}\left( n^{-d/(2m+2q+2)}\right)
 \end{eqnarray*}
 \end{thm}

\section{A simulation study}\label{simul}
 
\noindent  This section presents the results of simulations made in order to evaluate the performances of the proposed methods for slope estimation and  prediction in the model (\ref{rl1.1}).  We computed estimation and prediction errors from simulated spatial data in $\mathbb{Z}^2$.  Using the lexico-graphic order, we generated a sample  $\{(X_{\mathbf{i}_\ell},Y_{\mathbf{i}_\ell})\}_{1\leq  \ell\leq n^2}$ as follows:  we  consider the  $15$-th  first elements $B_1,\cdots,B_{15}$  of the B-splines basis. For $k=1,\cdots,15$, we generate a vector $(\xi_{\mathbf{i}_1,k},\cdots,\xi_{\mathbf{i}_{n^2},k})^T$ from a normal distribution $\mathcal{N}(0,\Sigma^1)$ in $\mathbb{R}^{n^2}$, where $\Sigma^1$ is the   $n^2\times n^2$ covariance matrix with general term  $\Sigma^1_{ij}=\exp(-3\|\mathbf{i}_i-\mathbf{i}_j\|_2)$. Further,  we  generate  a vector $(\Lambda_{\mathbf{i}_1} (t),\cdots,\Lambda_{\mathbf{i}_{n^2}}(t))^T$ from  a normal distribution $\mathcal{N}(0,\Sigma^2)$ in $\mathbb{R}^{n^2}$, where $\Sigma^2$ is the  $n^2\times n^2$ covariance matrix with general term  $\Sigma^2_{ij}=0.09$, and for $\ell=1,\cdots,n^2$ we take 
\[
X_{\mathbf{i}_\ell}(t)=\sum_{k=1}^{15}\xi_{\mathbf{i}_\ell,k}\,B_k(t)+\Lambda_{\mathbf{i}_\ell} (t).
\]
Considering $1001$ equispaced points in $[0,1]$, we compute each $Y_{\mathbf{i}_\ell}$  by approximating the integral in equation (\ref{rl1.1})  using   the rectangular method. That gives 
 \begin{eqnarray*}
 Y_{\mathbf{i}_\ell}&=\frac{1}{1000}&\sum^{1001}_{j=1}\beta(t_{j})X_{\mathbf{i}_\ell}(t_{j})+\epsilon_{\mathbf{i}_\ell}
 \end{eqnarray*}
 where $t_{j}=\frac{j-1}{1000}$, $j=1,\cdots,1001$,  the vector $(\epsilon_{\mathbf{i}_1},\cdots,\epsilon_{\mathbf{i}_{n^2}})^T$  is generated from a normal distribution $\mathcal{N}(0,\sigma^2_{\epsilon}\Sigma^1)$  with $\sigma^2_{\epsilon}$ controlled by the signal-to-noise ratio (snr) defined by \[
snr=\dfrac{\E[\langle\beta, X\rangle^2]}{\E[\langle\beta, X\rangle^2]+\sigma^2_{\epsilon}},
\]
and $\beta$ is a given function.  We considered two cases for the function $\beta$ given by:\

\bigskip

 \noindent\textbf{Case A :}  $\beta(t)=[\sin(2\pi t^3)]^3$;

 \noindent\textbf{Case B :} $\beta(t)=(0.4-t)^2$ .

\bigskip

\noindent The estimator $\widehat{\beta}$ of $\beta$ in model $(\ref{rl1.1})$ is computed by using the function "fregre.basis" of the $R$ fda package.  We assess performance of our methods  through the semi-norm $\|.\|_{\Gamma}$ defined in (\ref{s1}) for evaluating the estimation error between $\widehat{\beta}$ and $\beta$, and through the mean squared error (MSE) for evaluating the prediction error between the  prediction $\widehat{Y}_{\mathbf{i}_{0}}$ and the ''theoretical'' prediction $Y^{*}_{\mathbf{i}_{0}}$ at the  non-visited site $\mathbf{i}_0=(13.5,5)$. $X_{\mathbf{i}_0}$ is obtained by the ordinary krigging method, and $\widehat{Y}_{\mathbf{i}_{0}}$ and $Y^{*}_{\mathbf{i}_{0}}$ are obtained as defined in (\ref{s2}). We take $snr=5\%,10\%$ and $n=10,15,20,25$ over $100$ replications and we obtain the following tables.
    
   \bigskip

    \begin{center}
       \begin{tabular}{|c|c|c|c|c|c|}
       \hline
       snr(\%)&Case&$n^2=10^2$&$n^2=15^2$&$n^2=20^2$&$n^2=25^2$\\
       \cline{1-6}
        5 & A&0.073&0.0096&0.0091&0.0059\\\cline{2-6}
       & B&0.0440&0.0369&0.0176&0.0144\\\hline
       10& A&0.0085&0.0075&0.0069&0.0030\\\cline{2-6}
       & B&0.0352&0.0263&0.0202&0.0146\\\hline
       
       \end{tabular}

\bigskip

 \textbf{Table 1:} Estimation errors
    \end{center}

    \bigskip
    
        \begin{center}
           \begin{tabular}{|c|c|c|c|c|c|}
           \hline
           snr(\%)&Case&$n^2=10^2$&$n^2=15^2$&$n^2=20^2$&$n^2=25^2$\\
           \cline{1-6}
            5 & A&0.0008&0.0012&0.0011&0.0001\\\cline{2-6}
           & B&0.0097&0.0074&0.0024&0.0001\\\hline
           10& A&0.0008&0.0012&0.0011&0.0001\\\cline{2-6}
           & B&0.0055&0.0037&0.0021&0.0003\\\hline
           \end{tabular}
\bigskip

          \textbf{Table 2:}  Prediction errors at a non-visited site $\mathbf{i}_0=(13.5,5)$
        \end{center}

   \vspace{0.3cm}
   \noindent Table 1 and Table 2 present, respectively, the obtained estimation errors and prediction mean squared errors for different sample sizes and snr. The site $\mathbf{i}_0=(13.5,5)$ is beyond the grid of size $n^2=10^2$ whereas it is inside the grid of size $n^2=15^2$. We remark that when this point is inside the grid the prediction errors decrease as the sample size increases. Also, we see that estimation and prediction errors are small even when the sample size and the snr increase. 

\section*{Conclusion}
\noindent In this paper, we propose to study asymptotic properties of a smoothing splines estimator of slope function in a spatial functional linear regression model, where a scalar response is related to a square integrable spatial functional process. The originality of the proposed method is to consider spatially dependent data. 
 The main difficult is technical, especially in the proof of the prediction error because of the presence of the data spatial dependency. The prediction proposed in this work is available as well as for the points inside the grid than those beyond the grid compared to \cite{boukaetal18} where the prediction is only available for the points beyond the grid. One can then see the proposed methodology as a good alternative to \cite{Crambesetal09} when available data are spatially dependent. 

\section{Proofs}\label{section r5}

\subsection{A useful lemma}
Let
\begin{equation}\label{matrice}
\mathcal{M}=\left(\frac{1}{n^dp}\mathbf{X}^T\mathbf{X}+\rho\mathbf{A}_{m}\right)^{-1}
\left(\frac{1}{n^dp}\mathbf{X}^T\mathbf{X}\right);
\end{equation}
then we have:
\begin{lem}\label{lemme}
$ tr\left(\mathcal{M}^2\right)\leq\textrm{tr}\left(\mathcal{M}\right)$.
\end{lem}
\noindent\textit{Proof.} Since $\mathbf{A}_{m}$ is a symmetric nonnegative  matrix, its has a square root, denoted by  $\mathbf{A}_{m}^{1/2}$, that is also a symmetric nonnegative matrix. Denoting by  $\mathbf{A}_{m}^{-1/2}$ the inverse of   $\mathbf{A}_{m}^{1/2}$  and by $I_p$ the $p\times p$ identity matrix, we have:
\begin{eqnarray*}
\mathcal{M}
&=&\mathbf{A}_{m}^{-1/2}\left(\frac{1}{n^dp}\mathbf{A}_{m}^{-1/2}\mathbf{X}^{T}\mathbf{X}\mathbf{A}_{m}^{-1/2}+\rho I_p\right)^{-1}
\left(\frac{1}{n^dp}\mathbf{A}_{m}^{-1/2}\mathbf{X}^{T}\mathbf{X}\right).
\end{eqnarray*}
Then from the spectral decomposition
$
\frac{1}{n^dp}\mathbf{A}_{m}^{-1/2}\mathbf{X}^{T}\mathbf{X}\mathbf{A}_{m}^{-1/2}=\sum_{\ell=1}^p\mu_\ell\,u_\ell u_\ell^T,
$
where the $\mu_\ell$'s are the nonegative eigenvalues and $\left\{u_\ell\right\}_{1\leq\ell \leq p}$ is an orthonormal basis of $\mathbb{R}^p$ consisting of eigenvectors, it follows:
\begin{eqnarray*}
\mathcal{M}
&=&\sum_{\ell=1}^p\sum_{k=1}^p\frac{\mu_k}{\mu_\ell+\rho}\,\mathbf{A}_{m}^{-1/2}\,u_\ell u_\ell^Tu_k u_k^T\mathbf{A}_{m}^{1/2}
=\sum_{\ell=1}^p\frac{\mu_\ell}{\mu_\ell+\rho}\,\mathbf{A}_{m}^{-1/2}\,u_\ell u_\ell^T\mathbf{A}_{m}^{1/2}.
\end{eqnarray*}
Therefore, since $tr(\mathbf{A}_{m}^{-1/2}\,u_\ell u_\ell^T\mathbf{A}_{m}^{1/2})=tr(u_\ell^T\mathbf{A}_{m}^{1/2}\mathbf{A}_{m}^{-1/2}\,u_\ell )=tr(u_\ell^T\,u_\ell )=1$, we deduce that
$
tr(\mathcal{M})
=\sum_{\ell=1}^p\frac{\mu_\ell}{\mu_\ell+\rho}.
$
Finally,
\begin{eqnarray*}
tr\left(\mathcal{M}^2\right)
&=&tr\left(\sum_{\ell=1}^p\sum_{k=1}^p\left(\frac{\mu_\ell}{\mu_\ell+\rho}\right)\left(\frac{\mu_k}{\mu_k+\rho}\right)\,\mathbf{A}_{m}^{-1/2}\,u_\ell u_\ell^Tu_k u_k^T\mathbf{A}_{m}^{1/2}\right)\\
&=&\sum_{\ell=1}^p\left(\frac{\mu_\ell}{\mu_\ell+\rho}\right)^2\leq \sum_{\ell=1}^p\frac{\mu_\ell}{\mu_\ell+\rho}=tr(\mathcal{M}).
\end{eqnarray*}
\hfill $\Box$
\subsection{Proof of Theorem \ref{tr1}}
\noindent Putting
\begin{eqnarray*}
\Theta=\mathbf{X}\left(\frac{1}{n^dp}\mathbf{X}^T\mathbf{X}+\rho\mathbf{A}_{m}\right)^{-1}
\left(\frac{1}{n^dp}\mathbf{X}^T\mathbf{X}\right)\left(\frac{1}{n^dp}\mathbf{X}^T\mathbf{X}+\rho\mathbf{A}_{m}\right)^{-1}\mathbf{X}^T,
\end{eqnarray*}
we have
\begin{eqnarray}\label{erreurbeta}
\E_{\epsilon}\left(\left\|\widehat{\boldsymbol{\beta}}-\E_{\epsilon}(\widehat{\boldsymbol{\beta}})\right\|^{2}_{\Gamma_{n,p}}\right)
&=&\frac{1}{p}\E_\epsilon\left(\frac{1}{n^{2d}}\boldsymbol{\epsilon}^T\mathbf{X}\left(\frac{1}{n^dp}\mathbf{X}^T\mathbf{X}+\rho\mathbf{A}_{m}\right)^{-1}\right.\nonumber\\
&&\left.\left(\frac{1}{n^dp}\mathbf{X}^T\mathbf{X}\right)\left(\frac{1}{n^dp}\mathbf{X}^T\mathbf{X}+\rho\mathbf{A}_{m}\right)^{-1}\mathbf{X}^T\boldsymbol{\epsilon}\right)\nonumber\\
&=&\frac{1}{n^{2d}p}\left(\sum^{n^d}_{i=1}\Theta_{ii}\E\left(\tau_{i}^{2}\right)+\sum^{n^d}_{i=1}\sum_{\stackrel{j=1}{j\ne i}}^{n^d}\Theta_{ij}\E\left(\tau_{i}\tau_{j}\right)\right)
\end{eqnarray}
where $\tau_{i}=\epsilon_{\mathbf{i}_{i}}-\overline{\epsilon}$, with $\overline{\epsilon}=n^{-d}\sum_{j=1}^{n^d}\epsilon_{\mathbf{i}_{j}}$. Putting  $\sigma^{2}_{\epsilon}=\E\left(\epsilon^2_{\mathbf{i}_{i}}\right)$, we deduce from  $\tau^{2}_{i}=\epsilon^{2}_{\mathbf{i}_{i}}-2\epsilon_{\mathbf{i}_{i}}\overline{\epsilon}+\overline{\epsilon}^{2}$ and  the strict stationarity that
\begin{eqnarray}\label{etaui2}
\E\left(\tau^{2}_{i}\right)&=&\sigma^{2}_{\epsilon}-2\E\left[\epsilon_{\mathbf{i}_{i}}\left(\frac{1}{n^d}\sum^{n^d}_{j=1}\epsilon_{\mathbf{i}_{j}}\right)\right]+\E\left(\overline{\epsilon}^{2}\right)\nonumber\\
&=&\left(1-\frac{1}{n^d}\right)\sigma^{2}_{\epsilon}-\frac{2}{n^d}\sum_{\stackrel{j=1}{j\ne i}}^{n^d}\E\left(\epsilon_{\mathbf{i}_{i}}\epsilon_{\mathbf{i}_{j}}\right)+\frac{1}{n^{2d}}\sum_{k=1}^{n^d}\sum_{\stackrel{j=1}{j\ne k}}^{n^d}\E\left(\epsilon_{\mathbf{i}_{k}}\epsilon_{\mathbf{i}_{j}}\right)\nonumber\\
&\leq&\sigma^{2}_{\epsilon}+\frac{2}{n^d}\sum_{\stackrel{j=1}{j\ne i}}^{n^d}\left|\E\left(\epsilon_{\mathbf{i}_{i}}\epsilon_{\mathbf{i}_{j}}\right)\right|+\frac{1}{n^{2d}}\sum_{k=1}^{n^d}\sum_{\stackrel{j=1}{j\ne k}}^{n^d}\left|\E\left(\epsilon_{\mathbf{i}_{k}}\epsilon_{\mathbf{i}_{j}}\right)\right|.
\end{eqnarray}
Notice that, putting  $Q_n=\lfloor(\ln n)^{1/d}\rfloor$,  we have
\begin{eqnarray*}
\sum_{\stackrel{j=1}{j\ne i}}^{n^d}\left|\E\left(\epsilon_{\mathbf{i}_{i}}\epsilon_{\mathbf{i}_{j}}\right)\right|
=\sum_{\stackrel{j=1}{0<\delta(\{\mathbf{i}_{j}\},\{\mathbf{i}_{i}\})\leq Q_n}}^{n^d}\left|\E\left(\epsilon_{\mathbf{i}_{i}}\epsilon_{\mathbf{i}_{j}}\right)\right|+\sum_{\stackrel{j=1}{\delta(\{\mathbf{i}_{j}\},\{\mathbf{i}_{i}\})> Q_n}}^{n^d}\left|\E\left(\epsilon_{\mathbf{i}_{i}}\epsilon_{\mathbf{i}_{j}}\right)\right|.
\end{eqnarray*}
Then, using the Cauchy-Schwartz inequality as well as Lemma $2.1\ (ii)$ in \cite{tran},  we obtain,  under Assumption $\ref{ar2}$:
\begin{eqnarray*}
\sum_{\stackrel{j=1}{j\ne i}}^{n^d}\left|\E\left(\epsilon_{\mathbf{i}_{i}}\epsilon_{\mathbf{i}_{j}}\right)\right|
&\leq &\sigma_\epsilon^2\left(\sum_{\stackrel{j=1}{0<\delta(\{\mathbf{i}_{j}\},\{\mathbf{i}_{i}\})\leq Q_n}}^{n^d}1\right)\\
& +b_1&\sum_{\stackrel{j=1}{\delta(\{\mathbf{i}_{j}\},\{\mathbf{i}_{i}\})> Q_n}}^{n^d}\alpha_{1,\infty}\left(\delta(\{\mathbf{i}_{i}\},\{\mathbf{i}_{j}\})\right) \\
&\leq& \sigma^{2}_{\epsilon}\sum^{Q_n}_{t=1}t^{d-1}+b_{1}\sum^{\infty}_{k=Q_n+1}\,\,\sum_{k\leq t<k+1}\alpha_{1,\infty}\left(t\right) \\
&\leq&\sigma^{2}_{\epsilon}Q^{d}_n+b_{1}\sum^{\infty}_{t=Q_n+1}t^{d-1}\alpha_{1,\infty}\left(t\right) \\
&\leq& \sigma^{2}_{\epsilon}\ln( n)+b_{1}\sum^{\infty}_{t=1}t^{d-1-\theta},
\end{eqnarray*}
where $b_1$ is a positive constant. Since $\theta>d$, this finally gives:
\begin{equation}\label{majsom}
\sum_{\stackrel{j=1}{j\ne i}}^{n^d}\left|\E\left(\epsilon_{\mathbf{i}_{i}}\epsilon_{\mathbf{i}_{j}}\right)\right|
 \leq  \sigma^{2}_{\epsilon}\ln( n)+K_1,
\end{equation}
where $K_1$ is a positive constant. Therefore, from (\ref{etaui2}), it follows:
\begin{eqnarray}\label{etaui2_2}
\E\left(\tau^{2}_{i}\right)
&\leq&\sigma^{2}_{\epsilon}+\frac{3}{n^d}\left(\sigma^{2}_{\epsilon}\ln (n)+K\right).
\end{eqnarray}
Clearly,  $\sum^{n^d}_{i=1}\Theta_{ii}=tr(\Theta)=n^dp\,tr(\mathcal{M}^2)$, where $\mathcal{M}$ is defined in (\ref{matrice}). Then, from Lemma \ref{lemme}, and the proof of Theorem 1 in Crambes et al (2009) (see p. 55-56) that shows that $tr(\mathcal{M})\leq m+\rho^{-1/(2m+2q+1)}(2+C.C_{0})$ where $C$ and $C_0$ are positive constants, it follows:
\begin{eqnarray}\label{u}
\sum^{n^d}_{i=1}\Theta_{ii}&\leq&n^dp\left(m+\rho^{-1/(2m+2q+1)}(2+C.C_{0})\right).
\end{eqnarray}
Then, we deduce from (\ref{etaui2_2}) and (\ref{u}) that
\begin{eqnarray}\label{s1}
\frac{1}{n^{2d}p}\sum^{n^d}_{i=1}\Theta_{ii}\E\left(\tau_{i}^{2}\right)
\le\left(\frac{\sigma^{2}_{\epsilon}}{n^d}+\frac{c_{1}\ln ( n)}{n^{2d}}\right)\left(m+\rho^{-1/(2m+2q+1)}(2+C.C_{0})\right),
\end{eqnarray}
where $c_1$ is a positive constant. On the other hand,
\begin{eqnarray*}
\E\left(\tau_{i}\tau_j\right)&=&\E\left(\epsilon_{\mathbf{i}_{i}}\epsilon_{\mathbf{i}_{j}}-\epsilon_{\mathbf{i}_{i}}\overline{\epsilon}-\epsilon_{\mathbf{i}_{j}}\overline{\epsilon}+\overline{\epsilon}^2\right)
=\E\left(\epsilon_{\mathbf{i}_{i}}\epsilon_{\mathbf{i}_{j}}\right)-\frac{1}{n^d}\sigma^{2}_{\epsilon}\\
&-&\frac{1}{n^d}\sum_{\stackrel{k=1}{k\ne i}}^{n^d}\E\left(\epsilon_{\mathbf{i}_{i}}\epsilon_{\mathbf{i}_{k}}\right)-\frac{1}{n^d}\sum_{\stackrel{k=1}{k\ne j}}^{n^d}\E\left(\epsilon_{\mathbf{i}_{j}}\epsilon_{\mathbf{i}_{k}}\right)
+\frac{1}{n^{2d}}\sum_{k=1}^{n^d}\sum_{\stackrel{\ell=1}{\ell\ne k}}^{n^d}\E\left(\epsilon_{\mathbf{i}_{k}}\epsilon_{\mathbf{i}_{\ell}}\right).
\end{eqnarray*}
Then, using (\ref{majsom}), we obtain $\vert\E\left(\tau_{i}\tau_j\right)\vert\leq
\vert\E\left(\epsilon_{\mathbf{i}_{i}}\epsilon_{\mathbf{i}_{j}}\right)\vert +\frac{\sigma^{2}_{\epsilon}}{n^d}+\frac{3}{n^d}\left(\sigma^{2}_{\epsilon}\ln (n)+K_1\right)$, and 
\begin{eqnarray}\label{etauij}
\sum_{\stackrel{j=1}{j\ne i}}^{n^d}\vert\E\left(\tau_{i}\tau_j\right)\vert&\leq&
\sum_{\stackrel{j=1}{j\ne i}}^{n^d}\vert\E\left(\epsilon_{\mathbf{i}_{i}}\epsilon_{\mathbf{i}_{j}}\right)\vert +\sigma^{2}_{\epsilon}+3\left(\sigma^{2}_{\epsilon}\ln (n)+K_1\right)\nonumber\\
&\leq &\sigma^{2}_{\epsilon}+4\left(\sigma^{2}_{\epsilon}\ln (n)+K_1\right)
\leq K_2\ln(n),
\end{eqnarray}
where $K_2$ is a positive constant. Note that $\Theta=B^2$, where 
\[
B=(n^dp)^{-1}\mathbf{X}\left(\frac{1}{n^dp}\mathbf{X}^{T}\mathbf{X}+\rho\mathbf{A}_{m}\right)^{-1}
\mathbf{X}^{T};
\] 
then
\[
\vert\Theta_{ij}\vert=\vert\sum^{n^d}_{k=1}B_{ik}B_{kj}\vert\leq\frac{1}{2}\sum^{n^d}_{k=1}(B^{2}_{ik}+B^{2}_{kj})=\frac{1}{2}\left(\Theta_{ii}+\Theta_{jj}\right),
\]
and, putting $S=\frac{1}{n^{2d}p}\sum^{n^d}_{i=1}\sum_{\stackrel{j=1}{j\ne i}}^{n^d}\Theta_{ij}\E\left(\tau_{i}\tau_{j}\right)$, we deduce from this inequality and from (\ref{u}) and (\ref{etauij}) that
\begin{eqnarray*}
\vert S\vert&\leq &\frac{1}{n^{2d}p}\sum^{n^d}_{i=1}\sum_{\stackrel{j=1}{j\ne i}}^{n^d}\vert\Theta_{ij}\vert\,\vert\E\left(\tau_{i}\tau_j\right)\vert \leq\frac{1}{n^{2d}p}\sum^{n^d}_{i=1}\Theta_{ii}\,\sum_{\stackrel{j=1}{j\ne i}}^{n^d}\vert\E\left(\tau_{i}\tau_j\right)\vert\\
&\leq&\frac{K_2\ln(n)}{n^{2d}p}\sum^{n^d}_{i=1}\Theta_{ii}\leq \frac{K_2\ln(n)}{n^{d}}\left(m+\rho^{-1/(2m+2q+1)}(2+C.C_{0})\right).
\end{eqnarray*}
\hfill $\Box$\\

\subsection{Proof of  Theorem \ref{cr2}}

 \subsubsection{Lemma}
 
\begin{lem}\label{l4}
Under assumptions Theorem \ref{cr2}, we have
$$\|\widehat{\beta}-\beta\|^2=O\left(\frac{1}{\rho}\right)\ \ \text{a.s}.$$
\end{lem}

\subsubsection{Proofs}

{\bf Proof of Lemma \ref{l4} :}

\noindent We have
\begin{eqnarray*}
\|\widehat{\beta}-\beta\|^{2}\leq2\left(\|\widehat{\beta}\|^{2}+\|\beta\|^{2}\right)
\end{eqnarray*}
On the one hand, from assumption \ref{ar2}, we have  $\|\beta\|^{2}<C_{2}<\infty$. On the other hand, for $p$ large enough, we have
\begin{eqnarray*}
\|\widehat{\beta}\|^{2}&=&\left(\int^{1}_{0}\widehat{\beta}(t)^{2}dt-\frac{1}{p}\sum^{p}_{j=1}\widehat{\beta}(t_{j})^{2}\right)+\frac{1}{p}\widehat{\boldsymbol{\beta}}^{T}\widehat{\boldsymbol{\beta}}\le M_3+\frac{1}{p}\widehat{\boldsymbol{\beta}}^{T}\widehat{\boldsymbol{\beta}}
\end{eqnarray*}
Set $\mathbf{V}=(V_1-\overline{V},\cdots,V_{n^d}-\overline{V})^T$,  where $V_{\ell}=\int^{1}_{0}\beta(t)X_{\mathbf{i}_{\ell}}(t)dt-\frac{1}{p}\sum^{p}_{j=1}\beta(t_{j})X_{\mathbf{i}_{\ell}}(t_{j})$, $\ell=1,\cdots,n^d$. Then, by definition of $\widehat{\boldsymbol{\beta}}$, we have
\begin{eqnarray}
\frac{1}{p}\widehat{\boldsymbol{\beta}}^{T}\widehat{\boldsymbol{\beta}}&\leq&\frac{3}{p}\boldsymbol{\beta}^T\frac{1}{n^dp}\mathbf{X}^T\mathbf{X}\left(\frac{1}{n^dp}\mathbf{X}^T\mathbf{X}+\rho \mathbf{A}_m\right)^{-2}\frac{1}{n^dp}\mathbf{X}^T\mathbf{X}\boldsymbol{\beta}\nonumber\\
&+&\frac{3}{n^{2d}p}\mathbf{V}^T\mathbf{X}\left(\frac{1}{n^dp}\mathbf{X}^T\mathbf{X}+\rho \mathbf{A}_m\right)^{-2}\mathbf{X}^T\mathbf{V}\label{er17}\\
&+&\frac{3}{n^{2d}p}\boldsymbol{\epsilon}^T\mathbf{X}\left(\frac{1}{n^dp}\mathbf{X}^T\mathbf{X}+\rho \mathbf{A}_m\right)^{-2}\mathbf{X}^T\boldsymbol{\epsilon}\nonumber
\end{eqnarray}
The first and second term on the right-hand side of (\ref{er17}) are bounded as in \cite{Crambesetal09} (see p. 57), that is to say
$$\dfrac{3}{p}\boldsymbol{\beta}^T\dfrac{1}{n^dp}\mathbf{X}^T\mathbf{X}\left(\frac{1}{n^dp}\mathbf{X}^T\mathbf{X}+\rho \mathbf{A}_m\right)^{-2}\dfrac{1}{n^dp}\mathbf{X}^T\mathbf{X}\boldsymbol{\beta}=O(1)\ ,$$
$$\dfrac{3}{n^{2d}p}\mathbf{V}^T\mathbf{X}\left(\frac{1}{n^dp}\mathbf{X}^T\mathbf{X}+\rho \mathbf{A}_m\right)^{-2}\mathbf{X}^T\mathbf{V}=O\left(\dfrac{p^{-2\kappa}}{\rho}\right)$$
Set $W=\dfrac{1}{n^{d}p}\mathbf{X}\left(\frac{1}{n^dp}\mathbf{X}^T\mathbf{X}+\rho \mathbf{A}_m\right)^{-2}\mathbf{X}^T=BB^T$ where\\ $B=\dfrac{1}{\sqrt{n^{d}p}}\mathbf{X}\left(\frac{1}{n^dp}\mathbf{X}^T\mathbf{X}+\rho \mathbf{A}_m\right)^{-1}$. We have
\begin{eqnarray*}
N&:=&\frac{3}{n^{2d}p}\boldsymbol{\epsilon}^T\mathbf{X}\left(\frac{1}{n^dp}\mathbf{X}^T\mathbf{X}+\rho \mathbf{A}_m\right)^{-2}\mathbf{X}^T\boldsymbol{\epsilon}\\
&=&\frac{3}{n^{d}}\sum^{n^d}_{i,j=1}W_{ij}(\epsilon_{\mathbf{i}_i}-\overline{\epsilon})(\epsilon_{\mathbf{i}_j}-\overline{\epsilon})\\
&=&\frac{3}{n^{d}}\sum^{n^d}_{i=1}W_{ii}(\epsilon_{\mathbf{i}_i}-\overline{\epsilon})^2+\frac{3}{n^{d}}\sum^{n^d}_{{\stackrel{i,j=1}{i\ne j}}}W_{ij}(\epsilon_{\mathbf{i}_i}-\overline{\epsilon})(\epsilon_{\mathbf{i}_j}-\overline{\epsilon}):=N_1+N_2
\end{eqnarray*}
from Assumption \ref{ar22}, we have
$$|N_1|\le\frac{12M^2_1}{n^{d}}\sum^{n^d}_{i=1}W_{ii}=\frac{12M^2_1}{n^{d}}tr(W)\le \frac{12M^2_1}{n^{d}}tr[(\rho\mathbf{A}_m)^{-1}]=O\left(\frac{1}{n^d\rho}\right)\ \  \text{a.s}$$
and since
$$|W_{ij}|=\sum^{p}_{k=1}B_{ik}(B^T)_{kj}\le\frac{1}{2}\sum^{p}_{k=1}\left\{B_{ik}^2+\left[\left(B^T\right)_{kj}\right]^2\right\}=\frac{1}{2}(W_{ii}+W_{jj})$$
it follows that
$$|N_2|\le\frac{12M^2_1}{n^{d}}\sum^{n^d}_{{\stackrel{i,j=1}{i\ne j}}}|W_{ij}|\le 12M^2_1tr(W)=O\left(\frac{1}{\rho}\right)\ \  \text{a.s}$$
We then obtain the result of Lemma \ref{l4}.\hfill $\Box$

{\bf Proof of Theorem \ref{cr2} :}
\begin{eqnarray*}
B&:=&\E\left\{\E\left((\widehat{\beta}_{0}+\left\langle\widehat{\beta},X_{\mathbf{i}_{0}}\right\rangle-\beta_{0}-\left\langle\beta,X_{\mathbf{i}_{0}}\right\rangle)^{2}|\widehat{\beta}_{0}, \widehat{\beta}\right)\right\}\\
&=&\E\left\{\E\left[\left(\widehat{\beta}_{0}-\beta_{0}+\left\langle\widehat{\beta}-\beta,X_{\mathbf{i}_{0}}\right\rangle\right)^{2}|\widehat{\beta}_{0}, \widehat{\beta}\right]\right\}\\
&=&\E\left\{\E\left[\left(\left\langle\beta-\widehat{\beta},\overline{X}\right\rangle+\left\langle\widehat{\beta}-\beta,X_{\mathbf{i}_{0}}\right\rangle\right)^{2}|\widehat{\beta}_{0}, \widehat{\beta}\right]\right\}\\
&=&\E\left(\left\langle\widehat{\beta}-\beta,X_{\mathbf{i}_{0}}-\overline{X}\right\rangle^{2}\right)\\
&=&\E\left(\left(\left\langle\widehat{\beta}-\beta,X_{\mathbf{i}_{0}}-\E( X_{\mathbf{i}_{0}})\right\rangle+\left\langle\widehat{\beta}-\beta,\E(\overline{X})-\overline{X}\right\rangle\right)^{2}\right)\\
&\leq&2\left[\E\left(\left\langle\widehat{\beta}-\beta,X_{\mathbf{i}_{0}}-\E( X_{\mathbf{i}_{0}})\right\rangle^{2}\right)+\E\left(\left\langle\widehat{\beta}-\beta,\overline{X}-\E (\overline{X})\right\rangle^{2}\right)\right]:=B_1+B_2
\end{eqnarray*}
Since, from assumption \ref{ar4} and Lemma \ref{l4}, we have
$$\left|\left\langle\widehat{\beta}-\beta,X_{\mathbf{i}_{0}}-\E( X_{\mathbf{i}_{0}})\right\rangle\right|\le \|\widehat{\beta}-\beta\|\|X_{\mathbf{i}_{0}}-\E( X_{\mathbf{i}_{0}})\|=O(1/\sqrt{\rho})\ \ \text{a.s.}$$
and $\E\left(\left\langle X_{\mathbf{i}_{0}}-\E( X_{\mathbf{i}_{0}}),\zeta_j\right\rangle^4\right)\le4M^2_2\left\langle\Gamma\zeta_j,\zeta_j \right\rangle$,
it follows from Lemma $2.1\ (i)$ in \cite{tran}, Assumption \ref{ar4} and Lemma \ref{l4} that
\begin{eqnarray*}
B_1&:=&2\E\left(\left\langle\left\langle\widehat{\beta}-\beta,X_{\mathbf{i}_{0}}-\E( X_{\mathbf{i}_{0}})\right\rangle(X_{\mathbf{i}_{0}}-\E( X_{\mathbf{i}_{0}})),\widehat{\beta}-\beta\right\rangle\right)\\
&=&2\sum_{j\ge1}\E\left(\left\langle\widehat{\beta}-\beta,X_{\mathbf{i}_{0}}-\E( X_{\mathbf{i}_{0}})\right\rangle\left\langle X_{\mathbf{i}_{0}}-\E( X_{\mathbf{i}_{0}}),\zeta_j\right\rangle\left\langle\widehat{\beta}-\beta,\zeta_j\right\rangle\right)\\
&\le&2\sum_{j\ge1}\left\|\left\langle\widehat{\beta}-\beta,X_{\mathbf{i}_{0}}-\E( X_{\mathbf{i}_{0}})\right\rangle\left\langle X_{\mathbf{i}_{0}}-\E( X_{\mathbf{i}_{0}}),\zeta_j\right\rangle\right\|_4\left\|\left\langle\widehat{\beta}-\beta,\zeta_j\right\rangle\right\|_4\\
&&\times\left[\alpha_{1,\infty}(\delta(\{\mathbf{i}_0\},\{\mathbf{i}_1,...,\mathbf{i}_n\}))\right]^{1/2}\\
&&+2\sum_{j\ge1}\left|\E\left[\left\langle\widehat{\beta}-\beta,X_{\mathbf{i}_{0}}-\E( X_{\mathbf{i}_{0}})\right\rangle\left\langle X_{\mathbf{i}_{0}}-\E( X_{\mathbf{i}_{0}}),\zeta_j\right\rangle\right]\right|\left|\E\left[\left\langle\widehat{\beta}-\beta,\zeta_j\right\rangle\right]\right|\\
&\le&\frac{C}{\rho}\sum_{j\ge1}\left(\left\langle \Gamma\zeta_j,\zeta_j\right\rangle\right)^{1/4}\left[\alpha_{1,\infty}(\delta(\{\mathbf{i}_0\},\{\mathbf{i}_1,...,\mathbf{i}_n\}))\right]^{1/2}\\
&&+2\sum_{j\ge1}\left|\E\left[\left\langle\widehat{\beta}-\beta,X_{\mathbf{i}_{0}}-\E( X_{\mathbf{i}_{0}})\right\rangle\left\langle X_{\mathbf{i}_{0}}-\E( X_{\mathbf{i}_{0}}),\zeta_j\right\rangle\right]\right|\left|\E\left[\left\langle\widehat{\beta}-\beta,\zeta_j\right\rangle\right]\right|
\end{eqnarray*}
where $C$ is a positive constant. However, we have from Lemma $2.1\ (i)$ in \cite{tran} that
\begin{eqnarray*}
\Lambda&:=&\sum_{j\ge1}\left|\E\left[\left\langle\widehat{\beta}-\beta,X_{\mathbf{i}_{0}}-\E( X_{\mathbf{i}_{0}})\right\rangle\left\langle X_{\mathbf{i}_{0}}-\E( X_{\mathbf{i}_{0}}),\zeta_j\right\rangle\right]\right|\left|\E\left[\left\langle\widehat{\beta}-\beta,\zeta_j\right\rangle\right]\right|\\
&\le&\sum_{j\ge1}\sum_{\ell\ge1}\left\{\left|\E\left[\left\langle\widehat{\beta}-\beta,\zeta_{\ell}\right\rangle\left\langle X_{\mathbf{i}_{0}}-\E( X_{\mathbf{i}_{0}}),\zeta_{\ell}\right\rangle\left\langle X_{\mathbf{i}_{0}}-\E( X_{\mathbf{i}_{0}}),\zeta_j\right\rangle\right]\right|\right.\\
&&\left.\times \left|\E\left[\left\langle\widehat{\beta}-\beta,\zeta_j\right\rangle\right]\right|\right\}\\
&\le&\sum_{j\ge1}\sum_{\ell\ge1}\left\{\left\|\left\langle\widehat{\beta}-\beta,\zeta_{\ell}\right\rangle\right\|_4\left\|\left\langle X_{\mathbf{i}_{0}}-\E( X_{\mathbf{i}_{0}}),\zeta_{\ell}\right\rangle\left\langle X_{\mathbf{i}_{0}}-\E( X_{\mathbf{i}_{0}}),\zeta_j\right\rangle\right\|_4\right.\\
&&\left.\times\left|\E\left[\left\langle\widehat{\beta}-\beta,\zeta_j\right\rangle\right]\right|\right\}\left[\alpha_{1,\infty}(\delta(\{\mathbf{i}_0\},\{\mathbf{i}_1,...,\mathbf{i}_n\}))\right]^{1/2}\\
&&+\sum_{j\ge1}\sum_{\ell\ge1}\left\{\left|\E\left[\left\langle\widehat{\beta}-\beta,\zeta_{\ell}\right\rangle\right]\right|\left|\E\left[\left\langle X_{\mathbf{i}_{0}}-\E( X_{\mathbf{i}_{0}}),\zeta_{\ell}\right\rangle\left\langle X_{\mathbf{i}_{0}}-\E( X_{\mathbf{i}_{0}}),\zeta_j\right\rangle\right]\right|\right.\\
&&\left.\times\left|\E\left[\left\langle\widehat{\beta}-\beta,\zeta_j\right\rangle\right]\right|\right\}
\end{eqnarray*}
Since $\left\langle\Gamma\zeta_{\ell},\zeta_j \right\rangle=\lambda_{\ell}\left\langle\zeta_{\ell},\zeta_j \right\rangle=\lambda_{\ell}$ if $\ell=j$ and $0$ otherwise, it follows from Lemma \ref{l4} and Assumption \ref{ar4} that
$$\Lambda\le\frac{C_1}{\rho}\sum_{j\ge1}\left(\left\langle\Gamma\zeta_{\ell},\zeta_{\ell} \right\rangle\right)^{1/4}\left[\alpha_{1,\infty}(\delta(\{\mathbf{i}_0\},\{\mathbf{i}_1,...,\mathbf{i}_n\}))\right]^{1/2}+\left\|\E(\widehat{\beta}-\beta)\right\|^2_{\Gamma}$$
where $C_1$ is a positive constant. From assumption \ref{8} and Jensen inequality, we have
\begin{eqnarray*}
B_1\le\frac{C_2}{n^{d}\rho}\sum_{j\ge1}\lambda^{1/4}+2\left\|\E(\widehat{\beta}-\beta)\right\|^2_{\Gamma}\le\frac{C_3}{n^{d}\rho}+2\left\|\E(\widehat{\beta}-\beta)\right\|^2_{\Gamma}
\end{eqnarray*}
where $C_1$, $C_2$ and $C_3$ are positive constants. From Assumption \ref{ar4} and Lemma \ref{l4}, we have
\begin{eqnarray*}
B_2&:=&2\E\left(\left\langle\widehat{\beta}-\beta, \overline{X}-\E(\overline{X})\right\rangle^2\right)\\
&\le&\frac{K_1}{\rho}\E\left[\left\|\overline{X}-\E(\overline{X})\right\|^2\right]\\
&=&\frac{K_1}{n^{2d}\rho}\sum^{n^d}_{i=1}\E\left[\int^1_0\left( X_{\mathbf{i}_{i}}(u)-\E( X_{\mathbf{i}_{i}}(u))\right)^2 du\right]\\
&&+\frac{K_1}{n^{2d}\rho}\sum_{i\ne j}\int^1_0\E\left[\left( X_{\mathbf{i}_{i}}(u)-\E( X_{\mathbf{i}_{i}}(u))\right)\left(X_{\mathbf{i}_{j}}(u)-\E( X_{\mathbf{i}_{j}}(u))\right)\right]du\\
&\le&\dfrac{K_1}{n^d\rho}+\frac{K_2g(0)}{n^{2d}\rho}\sum^{n^d}_{i=1}\sum^{n^d}_{\stackrel{j=1}{j\ne i}}\varPsi(\delta(\{\mathbf{i}_i\},\{\mathbf{i}_{j}\}))\le\dfrac{K_1}{n^d\rho}+\frac{K_3}{n^d\rho}\sum^{\infty}_{t=1}t^{d-1}\varPsi(t)\leq \dfrac{K_4}{n^d\rho}
\end{eqnarray*}
where $K_1$, $K_2$, $K_3$ and $K_4$ are positive constants.  Then
$$B\leq \frac{K_5}{n^d\rho} + K_6\E\left[\|\widehat{\beta}-\beta\|^{2}_{\Gamma}\right]$$
where $K_7$ and $K_6$ are positive constants. Applying Corollary $\ref{cr1}$ with $2q>1$, $\rho\sim n^{-d(2m+2q+1)/(2m+2q+2)}$, we obtain the result of Theorem $\ref{cr2}$.
 
\hfill $\Box$

\end{document}